\def\IR{{\Bbb R}} 
 
\def\IS{{\Bbb S}}

\def\IC{\Bbb C} 
\def\ID{{\Bbb D}}

\def\oC{\hat{\IC}}

\documentclass[12pt]{article}

\usepackage[psamsfonts]{amssymb} 
\usepackage{amsfonts,amsmath}

\usepackage{epsfig,multicol}

\newtheorem{theorem}{Theorem}
\newtheorem{corollary}{Corollary}

\title{ Stream lines,  quasilines and holomorphic motions} 
\author{Gaven   J. Martin \thanks{Research supported by the Marsden Fund. \newline 
Mathematics Subject Classification Primary 30C62, 37F30, 30C75}} 
\date{} 
\begin{document}

\maketitle

\begin{abstract}  We give a new application of the theory of holomorphic motions to the study the distortion of level lines of harmonic functions and stream lines of ideal planar fluid flow. In various settings,  we show they are in fact quasilines - the quasiconformal images of the real line.  These methods also provide quite explicit global  estimates  on the geometry of these  curves.
\end{abstract}
    
\section{Introduction.} The theory of holomorphic motions, introduced by Ma{n}\'e-Sad-Sullivan \cite{MSS} and advanced by Slodkowski \cite{Slod},  has had a significant impact on the theory of quasiconformal mappings.  A reasonably thorough account of this is given in our book  \cite{AIM}.  In \cite{Martin, Martin2} we established some classical distortion theorems for quasiconformal mappings and used the theory to develop connections between Schottky's theorem and Teichm\"uller's theorem.  We also gave sharp estimates on the distortion of quasicircles which  in turn gave estimates for the distortion of extensions of analytic germs as studied in \cite{Jiang}.  Here we consider the geometry of stream lines for ideal fluid flow in a domain and establish bounds on their distortion in terms of a reference line.  These bounds come from an analysis of the geometry of the level lines of the hyperbolic metric and seem to be of independent interest.  When the reference line is known to be a quasiline - the   image of $\IR$ under a quasiconformal map of $\IC$ - which occurs for instance when there is some symmetry about,  it follows that all level lines are quasilines and it is possible to give explicit distortion estimates which contains global geometric information - such as bounded turning - for the curve,  see for instance  (\ref{ql}) below.  As such these estimates will have implications for parabolic linearisations.

\medskip

We first recall the two basic notions we  need here.  
\subsection{Quasiconformal mappings}
A homeomorphism $f:\Omega\to \IC$ defined on a domain $\Omega\subset \IC$ and in the Sobolev class $f\in W^{1,2}_{loc}(\Omega,\IC)$ of functions with locally square  integrable first derivatives is said to be  {\em quasiconformal} if there is a $1\leq K <\infty$ so that $f$ satisfies the {\em distortion inequality} 
\begin{equation}\label{de}
|Df(z)|^2 \leq K \; J(z,f), \hskip20pt\mbox{almost everywhere in $\Omega$}
\end{equation}
Here  $Df(z)$ is the Jacobian matrix and $J(z,f)$ its determinant.  If such a $K$ exists,  we will say that $f$ is  {\em $K$-quasiconformal.} The basic theory of quasiconformal mappings are  described in \cite{AIM}.  Quasiconformal mappings have the basic property of ``bounded distortion'' as they take infinitesimal circles to infinitesimal ellipses,  whereas conformal mappings have no distortion ($K=1$) as a consequence of the Cauchy-Riemann equations.  The distortion inequality (\ref{de}) actually implies the improved regularity $f\in W^{1,2K/(K+1)}_{loc}(\Omega)$,  \cite[Astala's theorem]{AIM}.

 \subsection{Holomorphic motions}
 
The theorem quoted below,  known as the extended $\lambda$-lemma and first proved by Slodkowski \cite{Slod},  is key in what follows.   The distortion estimate of $K$ in terms of the hyperbolic metric was observed by Bers \& Royden earlier and it is from this that we will be able to make our explicit distortion estimates below.  See \cite{GW} for a discussion.  A complete and accessible proof can be found in \cite[Chapter 12]{AIM}.  First,  the definition of a holomorphic motion.

\medskip
 
 Let $X\subset \oC=\IC\cup\{\infty\}$ be a set and $\ID$ the unit disk.   A holomorphic motion of $X$ is a map $\Phi:\ID\times X\to\oC$ such that
 \begin{itemize}
 \item For any fixed $a\in X$,  the map $\lambda\mapsto \Phi(\lambda,a)$ is holomorphic.
 \item For any fixed $\lambda\in \ID$,  the map $a\mapsto\Phi(\lambda,a)$ is an injection.
 \item $\Phi(0,a)=a$ for all $a\in X$.
 \end{itemize}
Note especially that there is no assumption regarding the measurability of $X$ or the continuity of $\Phi$ as a function of $a\in X$ or the two variables $(\lambda,a)\in\ID\times X$. 
 \begin{theorem}\label{holo} Let $\Phi:\ID\times X\to\oC$ be a holomorphic motion of $X$.  Then $\Phi$ has an extension to $\hat{\Phi}:\ID\times\oC$ which is a holomorphic motion of $\oC$ and for each $\lambda\in \ID$
\begin{equation} \label{Kest} \hat{\Phi}_\lambda = \hat{\Phi}(\lambda,\cdot):\oC\to\oC \hskip20pt\mbox{is $\frac{1+|\lambda|}{1-|\lambda|}$--quasiconformal.}
\end{equation}
 Moreover,  if $\rho_\ID$ denotes the hyperbolic metric (curvature = $-1$) of the unit disk,  then for $\lambda_1,\lambda_2\in \ID$ the map $\hat{\Phi}_{\lambda_1}^{-1} \circ \hat{\Phi}_{\lambda_2}$ is $K$--quasiconformal,  with $\log K= \rho_\ID(\lambda_1,\lambda_2)$.
 \end{theorem}
\noindent{\bf Remark:} We note that by using the Riemann mapping theorem here the parameter space $\ID$ can be replaced by any simply connected domain $\Omega$,  provided we replace $0$ by a point $\lambda_0\in \Omega$ and  assume that $\Phi(\lambda_0,a)=a$ for all $a\in X$.  Then the  estimate (\ref{Kest}) becomes
\begin{equation} \label{newKest}\hat{\Phi}_\lambda = \hat{\Phi}(\lambda,\cdot):\oC\to\oC \hskip20pt\mbox{is $e^{\rho_\Omega(\lambda_0,\lambda)}$--quasiconformal.} 
 \end{equation}
 
We will use this formulation when $\Omega$ is a strip.
    \section{Geometry of hyperbolic level lines.}

In what follows $\Omega\subset\IC$ will denote a Jordan domain, so that $\partial \Omega$ is a topological circle.  The reader will see that this is not necessary in much of what follows,  however it allows simplification in that we do not need to speak carefully of prime ends,  impressions and boundary values as every quasiconformal mapping between Jordan domains extends homeomorphically to the boundary.  In particular every Riemann map $\varphi:\ID\to\Omega$ extends homeomorphically to the boundary.  A hyperbolic line $\gamma \in \Omega$ is a complete hyperbolic geodesic.  In a Jordan domain $\gamma$ has two endpoints $\gamma_\pm \in \partial \Omega$.   For each such line there is a Riemann map $\varphi:\ID\to \Omega$ with $\varphi([-1,1])=\overline{\gamma}$ and $\varphi(-1)=\gamma_-$,  $\varphi(+1)=\gamma_+$.  For $c\geq 0$,  an arc $\alpha\subset\Omega$ is  a $c$-level line of $\gamma$ of the hyperbolic distance if for all $z\in \alpha$,
\begin{equation}
\rho_\Omega(z,\gamma)=c.
\end{equation}
Here $\rho_\Omega$ is the hyperbolic distance of $\Omega$.  Since,  by assumption,  $\Omega$ has a nice boundary the curve $\alpha$ extends to the boundary with endpoints $\alpha_\pm=\gamma_\pm$.  Our first theorem states that the level lines are the images of the line $\gamma$ under a self mapping of the domain with bounded distortion.

\begin{theorem} \label{ll} Let $\alpha$ be a $c$-level line of the hyperbolic line $\gamma\in \Omega$. Then there is a $K$-quasiconformal mapping $f:\Omega\to\Omega$ such that 
\begin{itemize}
\item $f(\gamma)=\alpha$,
\item $f|\partial \Omega = identity$,
\item $K \leq e^{c}$.
\end{itemize}
\end{theorem}
\noindent{\bf Proof.}  Let $\varphi:\ID\to\Omega$ with $\varphi((-1,1))= \gamma$.  There is a conformal mapping from the strip ${\cal S} = \{z\in\IC: |\Im m(z)| <\pi/2\}$ to $\ID$ so that the image of the real line is the interval $(-1,1)$.  Thus there is a conformal mapping $\psi:{\cal S} \to \Omega$ with $\psi(\IR)=\gamma$.  We define a holomorphic motion $\Phi(\lambda,a)$ of $\gamma$ and parameterised by $\lambda\in {\cal S}$ using the rule 
\begin{equation}
\Phi(\lambda,a) = \psi(\psi^{-1}(a)+\lambda)
\end{equation}
Then $\Phi(0,a)=identity$,  and this motion clearly depends holomorphically on $\lambda$ and is an injection for every $\lambda$.  The hyperbolic metric in ${\cal S}$ is $|dz|/\cos(y)$,  \cite[Example 7.9]{BeardonMinda},  thus the level lines (of $\IR$) of the hyperbolic distance are of the form $\{(x,t):x\in \IR\}$ and if we choose $t$ so that
\[ \int_{0}^{t} \frac{dy}{\cos(y)} = c, \]
then $\Phi(t,\gamma) = \alpha$.  Notice that for all $\lambda\in {\cal S}$, $\Phi(\lambda,\gamma) \subset \Omega$ and that $\Phi(\lambda,\gamma_{\pm}) = \gamma_\pm$.  Thus we can extend this function in the first instance by 
\begin{equation}
\Phi(\lambda,a) = \left\{ \begin{array}{ll} \psi(\psi^{-1}(a)+\lambda) & a\in \gamma  \\ a & a\in \partial\Omega \end{array} \right.
\end{equation}
and it is now clear that $\Phi$ defines a holomorphic motion of $\gamma\cup\partial\Omega$.   Theorem \ref{holo} now tells us that $\Phi$ is the restriction of a holomorphic motion of the Riemann sphere (note we have parametrised over ${\cal S}$ as per the remark following Theorem \ref{holo}).  At time $\lambda=t$ we obtain the quasiconformal mapping we seek and the distortion estimate follows from (\ref{newKest}) since $K\leq e^{\rho_{\cal S}(0,t)}$ and  $\rho_{\cal S}(0,t)=\rho_{\cal S}(\IR,\IR\pm t)=c$,  where $\rho_{\cal S}$ denotes the hyperbolic distance of ${\cal S}$.  \hfill $\Box$

\medskip

    { 
\scalebox{0.6}{\includegraphics[viewport= 10 600 720 750]{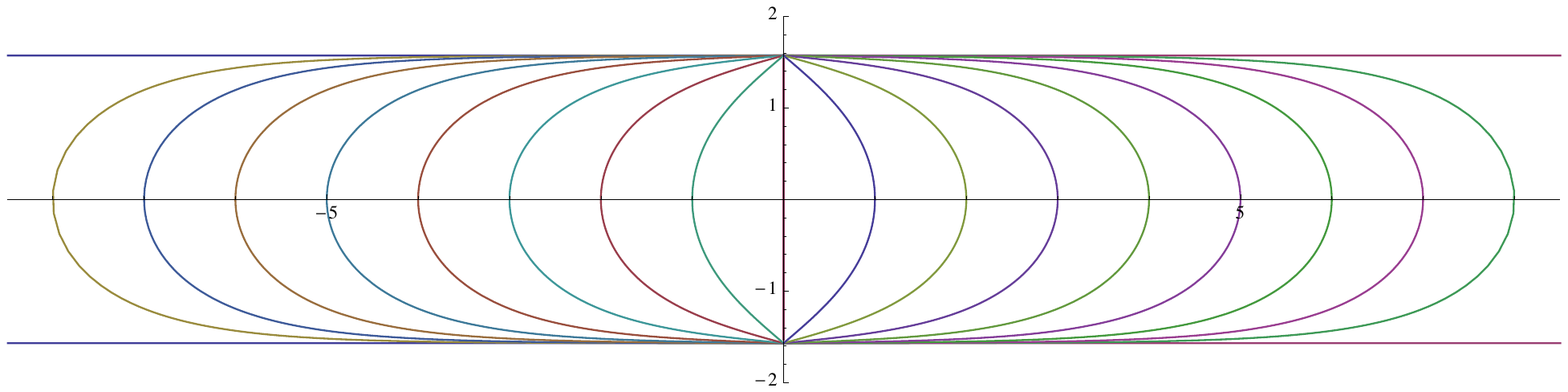}}
}
\noindent{\bf Fig. 1}  The integral level lines of the hyperbolic distance.  Theorem \ref{ll} gives a quasiconformal $f:\IC\to\IC$,  $f|\IC\setminus {\cal S}=identity$, $f$ is $e^{n}$-quasiconformal and takes the geodesic line   $[-\frac{i\pi}{2},\frac{i\pi}{2}]$ to the $n$-level line $\rho_{\cal S}(f([-\frac{i\pi}{2},\frac{i\pi}{2}]),[-\frac{i\pi}{2},\frac{i\pi}{2}])\equiv n$. 
\medskip

The next theorem concerns the level lines of the function associated with harmonic measure.

\begin{theorem} \label{thm3} Let $\Omega$ be a Jordan domain and $\alpha$ an open sub arc of $\partial \Omega$.  Let $h:\Omega\to\IR$ be the harmonic function with boundary values $h(z)=1$ for $z\in \alpha$ and $h(z)=0$ for $z\in \partial \Omega \setminus \bar \alpha$ .  Then for each $0<a\leq b < 1$, there is a $K$-quasiconformal mapping $f:\Omega\to\Omega$,  $f|\partial\Omega=identity$ and $f(\{h=a\}) = \{h=b\}$ with the distortion estimate
\begin{equation}\label{formula} K \leq  \frac{\tan \frac{b\pi}{2}}{\tan\frac{a\pi}{2}} 
\end{equation}
\end{theorem}
\noindent{\bf Remark:} We note that for  $b$ near $a=\frac{1}{2}$ we have
 \[ K_f \approx 1+ \pi\big|b-\frac{1}{2}\big| \]
 One might reasonably expect a quadratic term here though.
 
 \medskip
 
{\noindent \bf Proof.}  Let $\psi:{\cal S}\to\Omega$ be a conformal mapping with $\psi(\{y=\pi/2\}) = \alpha$ and $\phi(\{y=-\pi/2\}) = \partial\Omega\setminus \bar \alpha$.  A harmonic function which is $0$ on $\{y=-\frac{\pi}{2}\}$ and $1$ on $\{y=\frac{\pi}{2}\}$ is simply  $\tilde{h}(z)= \frac{y}{\pi}+\frac{1}{2}$, $z=x+iy$.  Its level lines are $\{\tilde{h}(z)= c\}=\{y=\frac{\pi}{2}(2c-1)\}$ for $0<c<1$.  Then $h(z)=\tilde{h}\circ \psi^{-1}$. The process is now as above,  we can holomorphically move these level lines keeping the boundary lines $\{y=\pm\frac{\pi}{2}\}$ and transfer this to $\Omega$ by $\psi$. The distortion is estimated by the exponential of the hyperbolic distance between the lines $\{h=a\}$ and $\{h=b\}$,  that is between the lines $\{y=\frac{\pi}{2}(2a-1)\}$ and $\{y =\frac{\pi}{2}(2b-1)\}$.  This is
\[ \rho = \Big|\int_{\frac{\pi}{2}(2a-1)}^{\frac{\pi}{2}(2b-1)} \frac{dt}{\cos(t)} \Big| = \log \frac{1+\tan(t/2)}{1-\tan(t/2)}\Big|_{t=\frac{\pi}{2}(2a-1)}^{t=\frac{\pi}{2}(2b-1)}  \]
and,  after a little manipulation,  we obtain the formula that we have given in (\ref{formula}).  \hfill $\Box$

\medskip

A {\em quasi-arc},  respectively {\em quasiline},  {\em quasicircle},  is the image of the line segment $(-1,1)$,  respectively  $\IR$,  $\IS=\partial \ID$, under a quasiconformal mapping $f:\IC\to\IC$.  Quasilines are simply quasicircles on the Riemann sphere which pass through $\infty$. When we know that the mapping $f$ is $K$-quasiconformal,  then we refer to $K$-quasiarcs, $K$-quasilines and $K$-quasicircles. There are many interesting geometric characterisations of these sets via criteria involving cross ratios and bounded turning initially discovered by Ahlfors,  see \cite{GH} for a comprehensive survey. For instance if $\alpha$ is a quasiline,  then there is a universal constant $C_\alpha$ such that if $z,w\in \alpha$ and $\alpha_{zw}$ denotes the finite subarc between them,  we have
\begin{equation}\label{ql}
{\rm diameter}(\alpha_{zw})\leq C_\alpha |z-w|
\end{equation}
This estimate is typically referred  to as ``bounded turning''.

\medskip

    { 
\scalebox{0.6}{\includegraphics[viewport= 10 550 720 790]{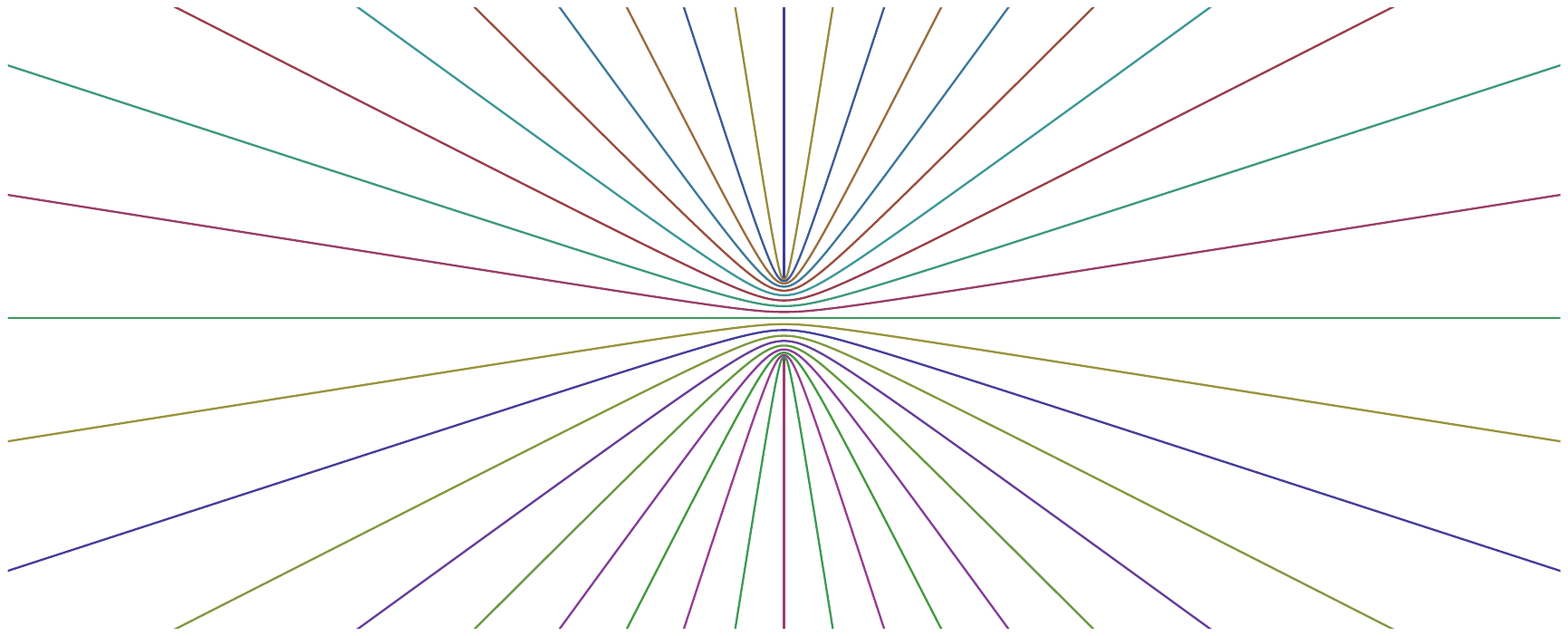}}
}\\
\noindent{\bf Fig. 2} Level lines of the harmonic function $h$ which is $1$ on $[\frac{i}{2},\infty)$ and $0$ on $[-\frac{i}{2},-\infty)$ and intervals of $\frac{1}{20}$.  $\IR=\{h=\frac{1}{2}$.  Theorem \ref{thm3} shows the level line $\{h=\frac{k}{20}\}$ to be a $\coth \frac{k\pi}{40}$-quasiline for $k\leq 10$ and $\tan \frac{k\pi}{40}$-quasiline for $10\leq k\leq 20$.

\medskip
\noindent{\bf Remark.}  Consideration of the level line $1-\delta$ as per Fig. 2 which is a $\tan\Big( \frac{(1-\delta)\pi}{2} \Big)$--quasiline with an asymptotic angle $\delta \pi/2$ with the imaginary axis at $\infty$, shows the bounded turning condition (\ref{ql}) cannot hold in general with a constant $C_\alpha < K$,  for a $K$-quasiline $\alpha$.

\medskip

We now have the next two corollaries.

\begin{corollary}  Let $\Omega$ and  $h$ be as in Theorem \ref{thm3}. If for some value of $a\in (0,1)$ the level line $\{h=a\}$ is a quasi-arc, then all level lines are quasi-arcs.  \end{corollary}

\begin{corollary}\label{cor2}  Let $\Omega$ be a Jordan domain symmetric across the real line,  $\alpha =\partial \Omega \cap \{\Im m(z)>0\}$, $\beta = \partial \Omega \cap \{\Im m(z)<0\}$ and $h$ the harmonic function with $h=1$ on $\alpha$,  $h=0$ on $\beta$.  Then for $b\geq \frac{1}{2}$ the level line 
\[ \gamma =\Big\{ \{h=b\} \cup \IR \setminus \Omega  \Big\} \]
is a quasi-line with distortion
\[ K \leq \tan \frac{b \pi}{2} \]
 \end{corollary}
 \noindent{\bf Proof.}  The hypotheses imply that $\{h=\frac{1}{2}\}$ is a segment of the real axis and thus a quasi-arc with distortion $K=1$.  The results follows once we extend the mapping $f$ given by Theorem \ref{thm3} by the identity outside $\Omega$ and note the resulting map is quasiconformal.  \hfill $\Box$

 \section{Ideal fluid flow.}  Knowing geometric information about the level lines of harmonic functions has many applications.  Here we give a couple which are simple and direct and concern ideal fluid flow in a channel and in particular two examples where computational results are known \cite{MM,WH}. In some ways the regularity results we derive (showing level lines are quasi lines) justify the computational results.
 
 \medskip
 
 A channel is the conformal image of the strip ${\cal S}$,  $\varphi:{\cal S}\to \IC$,  with the property that $\varphi(x+iy)\to \pm \infty$ as $x\to \pm \infty$.  A channel is not a Jordan domain,  but the reader can easily see that the above results apply without modification to this situation.

 Conformal invariance shows us that the stream lines of the fluid flow starting with a source at $-\infty$ and flowing to a sink at $+\infty$ are the level curves of a real valued nonconstant harmonic function $h$ which is constant on the boundary.  We can normalise so that these numbers are $h=+1$ on $\varphi(x+i\pi/2)$ and $h=-1$ on $\varphi(x-i\pi/2)$.  With this normalisation there is a central stream line $\alpha_0=\{h=0\}$.  We have the following two corollaries.

 \begin{corollary}  Let ${\cal C}$ be a channel and $\alpha$ a stream line for the flow of an incompressible fluid.  Then there is a quasiconformal map $f:\IC\to\IC$ with $f(\alpha_0)=\alpha$,  $f|\IC \setminus {\cal C}=identity$ and  
\[ K_f \leq e^{ d_{hyp}(\alpha_0,\alpha)} \]
where $ d_{hyp}$ is the hyperbolic metric of ${\cal C}$.
 \end{corollary}
 Of course $\alpha_0$ need not be a quasiline itself - the channel could be a regular neighbourhood of a smoothly embedded real line which is not quasiconformally equivalent to $\IR$ by a mapping of $\IC$.  The results says that every stream line is the bounded geometric image of the central line.
 Next,  if ${\cal C}$ is symmetric about the real line then $\alpha_0=\IR$ and we have the following.
\begin{corollary}\label{symcor}  Let ${\cal C}$ be a channel which is symmetric across the real line and $\alpha$ a stream line for the flow of an incompressible fluid.  Then $\alpha$ is a $K$-quasiline and   
\[ K  \leq e^{ d_{hyp}(\alpha,\IR)} \] 
 \end{corollary}
 
 We can also consider flow around an obstacle.  The following theorem has many obvious generalisations and we only present the simplest case.
 
 \begin{theorem}  Let $\Omega$ be a domain with $\overline{\Omega}\subset {\cal S}$.  Let $\alpha$ be a streamline homotopic to the either line $\{\Im m(z)= \pm\pi/2\}$ for ideal fluid flow in the obstructed channel ${\cal S}\setminus \Omega$.  Then $\alpha$ is a $K$-quasiline and 
 \[ K \leq e^{d_{hyp}(\alpha, \IR+i\pi/2)}\]
 where the metric is that of the domain   \[ D = \{  \Im m(z)=\pi/2\} \cup  \varphi(\{0<y<\pi/2\}) \cup (\varphi(\{0<y<\pi/2\})^{*} \} \]
 where $*$ denotes reflection in the line $\{\Im m(z)=\pi/2\}$.
 \end{theorem}
 \noindent{\bf Proof.}  The modulus of the ring ${\cal S}\setminus [-r,r]$ tends to infinity as $r\searrow 0$ and $0$ for $r=\infty$.  This modulus is continuous and so the intermediate value theorem gives us an $r$ so that ${\rm Mod}({\cal S}\setminus [-r,r])={\rm Mod}({\cal S}\setminus \Omega)$.  For this $r$ there is a conformal mapping $\varphi : {\cal S}\setminus [-r,r] \to {\cal S}\setminus \Omega$.  The stream lines for ideal fluid flow in  ${\cal S}\setminus [-r,r]$ are simply the lines $\{\IR+iy: 0< |y| \leq \pi/2 \}$. The images of these lines under $\varphi$ are the stream lines for flow in ${\cal S}\setminus \Omega$.  We can also use the Carath\'eodory/Schwarz extension/reflection principle to extend $\varphi$ to a map $\tilde{\varphi} : \{z:0< y < \pi \}\to \IC$,  $\tilde{\varphi}(\IR+i\pi/2)=\IR+i\pi/2$.  It follows from Corollary \ref{symcor} that for $0<y<\pi/2$,   the stream line $\tilde{\varphi}(\IR+iy)$ is a quasiline with distortion $e^{d_{hyp}(\tilde{\varphi}(\IR+iy), \IR+i\pi/2)}$ - with the metric here being that of $\tilde{\varphi}(\{0<y<\pi\})$. \hfill $\Box$
  
    { 
\scalebox{0.45}{\includegraphics[viewport= 20 350 420 580]{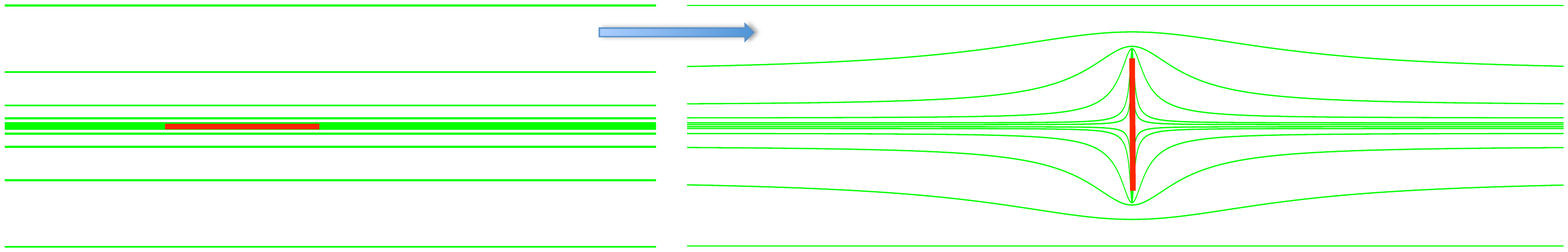}}
}
\\
\noindent{\bf Fig. 3}  The integral flow lines around the segment $[-i,i]$ in the strip $\{|\Im m(z)|<\frac{\pi}{2}\}$.

\medskip

It is not difficult to further refine this estimate upon consideration of the harmonic function involved following the arguments given for Corollary \ref{cor2}.  What is remarkable here is that the global geometric estimates one achieves on the stream lines do not depend on the complexity of the object.  Indeed the fact that the hyperbolic metric increases under inclusion implies that we can use a slightly larger smooth approximation to the object to get the estimates we require.  So the boundary being highly irregular (say Hausdorff dimension $>1$)  does not matter for estimating the distortion unless the stream line comes very close to the boundary.  Further,  the same argument (with the same estimates) works for flow around multiple objects,  though there are issues are there is more than one module.   In the simplest case,  the bounds on the distortion on the stream lines applies for flow around $\varphi({\cal S}\setminus  \cup I_i)\subset {\cal S}$,  where $\{I_i\}$ is {\em any} disjoint collection of closed intervals of $\IR$,  $\varphi$ is conformal with $\varphi(\{\Im m(z)=\pm\frac{\pi}{2}\}) = \{\Im m(z)=\pm\frac{\pi}{2}\}$

\bigskip

\noindent G.  Martin -  Massey University,  Auckland, NZ,    g.j.martin@massey.ac.nz

\end{document}